\documentclass[12pt]{amsart}
\usepackage{amssymb, amsmath, amsthm} 
\usepackage{url}      
\usepackage{mathtools} 
\usepackage[dvipsnames]{xcolor}
\usepackage{tikz-cd}
\usepackage[T1]{fontenc}
\usepackage[utf8]{inputenc}
\usepackage{mathrsfs}
\usepackage{stmaryrd}
\usepackage{eucal}
\usepackage[all]{xy}
\usepackage[margin=1in]{geometry} 
\usepackage{hyperref}
\usepackage[shortlabels]{enumitem}
\hypersetup{bookmarksdepth=2}
\hypersetup{colorlinks=true}
\hypersetup{linkcolor=blue}
\hypersetup{citecolor=blue}
\hypersetup{urlcolor=blue}
\usepackage[backend=bibtex,style=alphabetic]{biblatex}
\AtBeginBibliography{\raggedright}
\allowdisplaybreaks
\theoremstyle{plain}
\newtheorem{theorem}{Theorem}[section] 
\newtheorem{lemma}[theorem]{Lemma}

\newtheorem*{mainrestatement}{Theorem~\ref{thm:main}}

\theoremstyle{definition}
\newtheorem{definition}[theorem]{Definition}

\theoremstyle{remark}
\newtheorem{remark}[theorem]{Remark}
\newtheorem*{acknowledgement}{Acknowledgements}
\newtheorem*{ai}{AI Usage}

\renewcommand{\a}{\alpha}

\renewcommand{\r}{\rho}

\newcommand{\G}{\Gamma}

\newcommand{\gp}{{\mathfrak{p}}}

\newcommand{\Kbar}{{\overline K}}

\newcommand{\Ocal}{{\mathcal O}}

\newcommand{\Ucal}{{\mathcal U}}

\newcommand{\CC}{\mathbb{C}}
\newcommand{\FF}{\mathbb{F}}

\newcommand{\PP}{\mathbb{P}}

\newcommand{\ZZ}{\mathbb{Z}}

\newcommand{\characteristic}{\operatorname{char}}

\newcommand{\GL}{\operatorname{GL}}

\newcommand{\grand}{{\textup{grand}}} 
\newcommand{\grano}{\mathcal{O}^{\grand}}

\newcommand{\id}{\operatorname{id}}

\newcommand{\Mat}{\operatorname{Mat}}

\title{Orbit Transversality for Projective Space}

\author{Kaiwen Lu}
\address{Department of Mathematics, Brown University, Providence, RI, USA}
\email{kaiwen\_lu@brown.edu}

\begin{document}

\begin{abstract}
Let \(f\colon \PP^n\to \PP^n\) be a dominant endomorphism defined over a number field \(K\).
We show that there is a finite extension \(L/K\) such that every proper closed subset of \(\PP^n\) misses at least one grand \((f,L)\)-orbit.
Equivalently, \((\PP^n,f,K)\) is strongly dense orbit transversal (S-DOT): every set containing exactly one representative of each grand \((f,L)\)-orbit is Zariski dense in \(\PP^n\).
\end{abstract}

\maketitle

\section{Introduction}
The classical Mordell--Lang conjecture, proved by Faltings~\cite{faltings} for abelian varieties in characteristic zero, describes the intersection of an algebraic subvariety with a finitely generated subgroup as a finite union of cosets of subgroups.
Replacing the finitely generated subgroup by the forward orbit of a point under a self-morphism motivates the Dynamical Mordell--Lang conjecture, which does not require a group structure. 
For a quasiprojective variety \(X\) over a field \(K\) of characteristic zero and a \(K\)-morphism \(f\colon X\to X\), it predicts that, for any point \(x\in X(K)\) and closed subvariety \(V\subseteq X\), the return-time set
\begin{equation*}
    \{n\in\ZZ_{\ge 0}\mid f^n(x)\in V(K)\}
\end{equation*}
is a finite union of arithmetic progressions, including singleton sets~\cite{bellghiocatucker}.

Dynamical Mordell--Lang concerns how one forward orbit intersects a closed subvariety. 
We consider a complementary question, varying the point and asking whether an entire grand orbit can avoid a given proper closed subset:
\begin{quote}
\itshape
Can one choose a finite extension of the ground field over which every proper Zariski-closed subset misses at least one rational grand orbit?
\end{quote}
We make this question precise using grand orbits and transversals.

Let \(X\) be a variety over a field \(K\) and let \(f\colon X\to X\) be a morphism defined over \(K\).
\begin{definition}
    The \emph{grand \((f,K)\)-orbit} of a point \(P\in X(K)\) is the set of points whose orbits eventually merge with the orbit of \(P\):
    \[\grano_{f,K}(P):=\{Q\in X(K)\mid  f^n(Q)=f^m(P) \text{ for some \(n,m \ge 0\)}\}.\]
\end{definition}
\begin{definition}
An \emph{\((f,K)\)-transversal} is a subset of \(X(K)\) containing exactly one representative of each grand \((f,K)\)-orbit.
We say that \((X,f,K)\) is \emph{strongly dense orbit transversal} (S-DOT) if for some finite extension \(L/K\), \textbf{every} \((f,L)\)-transversal of \(X\) is Zariski dense.
\end{definition}

The notion of S-DOT can be formulated with respect to any topology on the rational points. Throughout this paper, S-DOT refers to the Zariski topology.

A proper closed subset meets every grand \((f,K)\)-orbit exactly when it contains an \((f,K)\)-transversal. Thus \((X,f,K)\) is S-DOT precisely when, over one finite extension \(L/K\), no proper closed subset meets every grand \((f,L)\)-orbit.

Our main result gives an affirmative answer to the orbit avoidance question for dominant endomorphisms of projective space over number fields.
\begin{theorem}\label{thm:main}
    Let \(f\colon \PP^n\to \PP^n\) be a dominant endomorphism over a number field \(K\). Then there is a finite extension \(L/K\) such that for every proper closed subset \(V\subsetneq \PP^n\), there is a point \(P\in \PP^n(L)\) such that 
    \[\grano_{f,L}(P)\cap V(L)=\varnothing.\]
    In other words, \((\PP^n,f,K)\) is S-DOT.
\end{theorem}

This question was first posed by Pasten and Silverman in their study of orbit propagation~\cite{pastensilverman}. 
The discussion above gives a different motivation for the same problem. Their work studies how the existence of one Zariski dense forward orbit is reflected in the distribution of grand orbits. 
For projective space, under the dense orbit hypothesis, their results imply S-DOT for linear endomorphisms and the existence of a Zariski dense transversal for endomorphisms of degree at least~\(2\)~\cite[Theorem~4.1]{pastensilverman}.

They also asked whether potential density of rational points, meaning Zariski density over some finite extension, suffices for S-DOT~\cite[Remark~1.9]{pastensilverman}. 
Since \(\PP^n(K)\) is Zariski dense, Theorem~\ref{thm:main} answers this question affirmatively for dominant endomorphisms of projective space. In degree at least~\(2\), it strengthens their conclusion from the existence of a dense transversal to the density of every transversal. 
We established an analogous result in our earlier preprint~\cite{sdot}, proving S-DOT for dominant endomorphisms of arbitrary abelian varieties without assuming the existence of a Zariski dense orbit.

\begin{remark}
    Over a finite extension \(L/K\), a proper closed subset can also be viewed as capturing the rational dynamics if it meets every grand \((f,L)\)-orbit. Equivalently, it sweeps out all of \(X(L)\) by taking forward images of its \(L\)-rational points and then iterated preimages in \(X(L)\). If \(f\) is not dominant, the proper closed subset \(\overline{f(X)}\) does meet every grand \((f,L)\)-orbit, so dominance is necessary for S-DOT.
\end{remark}

\begin{ai}
    ChatGPT 6 Astra was used to audit the drafts. In addition to typographical and grammatical corrections, Astra suggested expressing \(\G\) as in equation~\eqref{gammaform} instead of using a multivariable Taylor expansion, which simplified the proof.
\end{ai}
\begin{acknowledgement}
    I would like to thank my advisor Joseph Silverman for many helpful discussions and comments.
\end{acknowledgement}

\section{Structure of the proof}
The structure of the proof of the main theorem is as follows:
\begin{enumerate}
    \item When \(\deg(f)\ge 2\), we use Lemma~\ref{lem:algperpt} to find an algebraic repelling periodic point \(\a\), and Lemma~\ref{lem:iterinv} to replace \(f\) by an iterate so that \(\a\) becomes a repelling fixed point. Here \(\deg(f)\) denotes the algebraic degree.
    \item At a suitable prime with uniformizer \(\pi\), we use Lemma~\ref{lem:hypredcount} to bound the number of reductions of rational points on a hypersurface modulo \(\pi^m\), and Lemma~\ref{lem:finitedyn} to bound the order of the permutation induced by \(f\) modulo \(\pi^m\) on a ball around \(\a\).
    \item The map \(f\) restricts to a self-map on a \(\pi\)-adic ball around \(\a\). To control grand orbits meeting this ball, we give a uniform bound on the first entrance time of any orbit that eventually enters it.
    \item We then combine the hypersurface and cycle-length bounds in a counting argument to find an avoiding grand orbit through a point in the ball. For maps of degree \(1\), we apply the counting argument to the permutation induced on the entire projective space modulo \(\pi^m\).
\end{enumerate}

\section{Periodic points and reduction estimates}
We first justify passing to an iterate and choosing an algebraic repelling periodic point. The remaining lemmas bound reductions of hypersurfaces and cycle lengths.

\begin{lemma}\label{lem:iterinv}
    Let \(K\) be a field, \(X/K\) an irreducible variety, and
    \(f\colon X\to X\) an endomorphism defined over \(K\).
    Then for all \(n\ge 1\), every \((f,K)\)-transversal is dense in \(X\)
    if and only if every \((f^n,K)\)-transversal is dense in \(X\).

    If \(X\) is geometrically irreducible, this implies that \((X,f,K)\)
    is S-DOT if and only if \((X,f^n,K)\) is S-DOT.
\end{lemma}
\begin{proof}
    Fix \(n\ge 1\). We prove both implications by contraposition.
    Let \(S'\subseteq X(K)\) be an \((f^n,K)\)-transversal such that
    \(\overline{S'}\) is a proper closed subset of \(X\).
    The set \(S'\) contains an \((f,K)\)-transversal, because every grand
    \((f,K)\)-orbit is a union of grand \((f^n,K)\)-orbits, each of which
    has a representative in \(S'\).
    Therefore we may take a subset \(S\subseteq S'\) that is an
    \((f,K)\)-transversal, and \(S\) is not dense in \(X\).

    Conversely, let \(S\subseteq X(K)\) be an \((f,K)\)-transversal such that
    \(\overline{S}\) is a proper closed subset of \(X\).
    We claim that \(\bigcup_{i=0}^{n-1}f^i(S)\) contains an
    \((f^n,K)\)-transversal.
    Indeed, for every \(x\in X(K)\), there are \(s\in S\) and integers
    \(a,b\ge 0\) such that
    \[
        f^a(x)=f^b(s).
    \]
    Choose \(r\in\{0,\ldots,n-1\}\) such that \(n\) divides \(a+r\),
    and write \(b+r=n\ell+i\), where \(\ell\ge 0\) and \(0\le i<n\).
    Then
    \[
        f^{a+r}(x)=f^{b+r}(s)=f^{n\ell}(f^i(s)).
    \]
    Since \(n\) divides \(a+r\), this shows that \(x\) and \(f^i(s)\)
    are in the same grand \((f^n,K)\)-orbit.
    Thus every grand \((f^n,K)\)-orbit meets
    \(\bigcup_{i=0}^{n-1}f^i(S)\), and we may take an
    \((f^n,K)\)-transversal
    \[
        S'\subseteq \bigcup_{i=0}^{n-1}f^i(S).
    \]

    By continuity of the maps \(f^i\), we have
    \[
        \overline{\bigcup_{i=0}^{n-1}f^i(S)}
        =\bigcup_{i=0}^{n-1}\overline{f^i(S)}
        =\bigcup_{i=0}^{n-1}\overline{f^i(\overline{S})}.
    \]
    Each set \(\overline{f^i(\overline{S})}\) is a proper closed subset
    of \(X\), since
    \[
        \dim\overline{f^i(\overline{S})}
        \le \dim\overline{S}<\dim X.
    \]
    Since \(X\) is irreducible, the finite union of these proper closed
    subsets cannot equal \(X\).
    Then \(\overline{S'}\ne X\), and \(S'\) is not dense.

    Finally, if \(X\) is geometrically irreducible, then \(X_L\) is
    irreducible for every finite extension \(L/K\).
    The same argument therefore gives the equivalence over every
    such \(L\), proving the assertion about S-DOT.
\end{proof}

\begin{lemma}\label{lem:algperpt}
    Let \(K\) be a field with a fixed embedding \(\Kbar\hookrightarrow\CC\),
    and let \(f\colon \PP^n\to\PP^n\) be an endomorphism defined over \(K\)
    of algebraic degree \(d\ge 2\), where \(n\ge 1\).
    Then there exists a point \(P\in\PP^n(\Kbar)\) such that \(P\) is
    a repelling periodic point with respect to the chosen embedding.
    That is, there exists \(r\ge 1\) such that \(f^r(P)=P\) and all
    eigenvalues of \(d(f^r)_P\) have complex absolute value greater than \(1\).
\end{lemma}
\begin{proof}
    Using the chosen embedding, consider the map
    \(f\colon \PP^n(\CC)\to\PP^n(\CC)\).
    By \cite[Theorem~2]{briendduval}, the repelling periodic points of a holomorphic endomorphism of \(\PP^n(\CC)\) of algebraic degree at least \(2\) equidistribute with respect to the equilibrium measure, and in particular exist.
    Let \(P\in\PP^n(\CC)\) be a repelling periodic point with exact period \(r\).
    Then \(P\) is a fixed point of \(f^r\).
    It remains to show that \(P\) is algebraic over \(K\).
    We do this by showing that the fixed-point scheme of \(f^r\) is finite over \(K\).

    Let \(Y\subseteq\PP^n\) be the fixed-point scheme of \(f^r\).
    Since \(f\) is defined over \(K\), so is \(Y\).
    Suppose for contradiction that \(Y_{\Kbar}\) has positive dimension.
    Then it contains an integral projective curve \(C\).
    Let \(i\colon C\hookrightarrow\PP^n_{\Kbar}\) be the inclusion.
    Since \(C\) is contained in the fixed-point scheme of \(f^r\),
    we have \(f^r\circ i=i\).
    On the other hand, the algebraic degree of \(f^r\) is \(d^r\), so
    \[
        (f^r)^*\mathcal{O}_{\PP^n_{\Kbar}}(1)
        \cong \mathcal{O}_{\PP^n_{\Kbar}}(d^r).
    \]
    Pulling back to \(C\) gives
    \[
        \mathcal{O}_C(1)
        =i^*\mathcal{O}_{\PP^n_{\Kbar}}(1)
        \cong i^*(f^r)^*\mathcal{O}_{\PP^n_{\Kbar}}(1)
        \cong \mathcal{O}_C(d^r).
    \]
    Taking degrees on \(C\) gives
    \[
        \deg\mathcal{O}_C(1)
        =d^r\deg\mathcal{O}_C(1),
    \]
    which is impossible because \(\deg\mathcal{O}_C(1)>0\) and \(d^r>1\).

    Thus \(Y_{\Kbar}\), and hence \(Y\), is zero-dimensional.
    Since \(Y\) is a closed subscheme of \(\PP^n_K\), it is projective
    over \(K\), and therefore finite over \(K\).
    Consequently, every complex point of \(Y\) is defined over a finite
    extension of \(K\), so \(P\in\PP^n(\Kbar)\).
\end{proof}

\begin{lemma}\label{lem:hypredcount}
    Let \(A\) be a discrete valuation ring with a uniformizer \(\pi\), fraction field \(K\), and a finite residue field \(k\) with \(\#k=q\) and \(\characteristic(k)=p\), and let \(F\in K[x_1,\ldots,x_n]\) be a nonzero polynomial of total degree at most \(d\), where \(d\ge 1\).
    For each integer \(m\ge 1\), consider the set of exact zeros and its
    image under the reduction map \(\r_m\colon A^n\to (A/\pi^m A)^n\):
    \[
        Z:=\{x\in A^n\mid F(x)=0\}, \quad Z_m:=\r_m(Z).
    \]
    Then we find
    \[
        \#Z_m \le (dq^{n-1})^m.
    \]
\end{lemma}
\begin{proof}
    We use the valuation on \(K\) normalized so that \(\pi\) has valuation \(1\), and proceed by induction on \(m\).
    When \(m=1\), we have \(A/\pi A=k\).
    Let \(v\) be the minimum of the valuations of the nonzero coefficients
    of \(F\), so
    \[
        F'=\pi^{-v}F\in A[x_1,\ldots,x_n].
    \]
    This polynomial has at least one coefficient that is a unit in \(A\),
    and it has the same zero set as \(F\).
    By \cite[Theorem~6.13]{finitefields}, the nonzero reduction
    \(\overline{F'}\in k[x_1,\ldots,x_n]\) of \(F'\) has at most \(dq^{n-1}\)
    zeros in \(k^n\).
    Since \(Z_1\) is contained in the zero set of \(\overline{F'}\) in \(k^n\),
    we find \(\#Z_1\le dq^{n-1}\).

    Now let \(m\ge 2\), and assume that the bound is true for \(m-1\).
    It suffices to show that for each point in \(Z_{m-1}\), at most
    \(dq^{n-1}\) of its preimages under the quotient map
    \[
        (A/\pi^m A)^n\to (A/\pi^{m-1}A)^n
    \]
    are images of actual zeros of \(F\).

    Let \(a+\pi^{m-1}A^n\) be a class in \(Z_{m-1}\) with \(a\in A^n\),
    and consider \(G'\in K[y_1,\ldots,y_n]\) defined by
    \[
        G'(y)=F(a+\pi^{m-1}y).
    \]
    The polynomial \(G'\) is nonzero, because substitution by
    \(a+\pi^{m-1}y\) is an invertible affine change of variables over \(K\),
    and its degree is at most \(d\).
    Take \(\nu\) to be the minimum of the valuations of the nonzero
    coefficients of \(G'\), so
    \[
        G=\pi^{-\nu}G'\in A[y_1,\ldots,y_n]
    \]
    and at least one coefficient of \(G\) is a unit in \(A\).
    Thus \(G\) has the same zero set as \(G'\), and its reduction
    \(\overline{G}\in k[y_1,\ldots,y_n]\) is a nonzero polynomial of degree at most \(d\).

    Note that the elements in \((A/\pi^m A)^n\) that reduce to
    \(a+\pi^{m-1}A^n\) are precisely those of the form
    \(\r_m(a+\pi^{m-1}b)\) with \(b\in A^n\).
    These elements are parameterized bijectively by \(\overline{b}\in k^n\),
    since two choices of \(b\) give the same element if and only if they
    agree modulo \(\pi A^n\).
    If such an element is the image of an actual zero \(x\in Z\) of \(F\),
    then we may replace \(b\) by \((x-a)/\pi^{m-1}\in A^n\) without changing
    \(\overline{b}\).
    With this choice, we have
    \[
        G(b)=\pi^{-\nu}G'(b)
        =\pi^{-\nu}F(a+\pi^{m-1}b)=0,
    \]
    and reducing modulo \(\pi\) gives \(\overline{G}(\overline{b})=0\).
    By \cite[Theorem~6.13]{finitefields}, there are at most \(dq^{n-1}\)
    such \(\overline{b}\in k^n\), since \(\overline{G}\) is nonzero and
    has degree at most \(d\).
    Thus each point in \(Z_{m-1}\) has at most \(dq^{n-1}\) preimages in
    \(Z_m\), so
    \[
        \#Z_m
        \le (dq^{n-1})\#Z_{m-1}
        \le (dq^{n-1})(dq^{n-1})^{m-1}
        =(dq^{n-1})^m,
    \]
    finishing the induction.
\end{proof}

\begin{remark}
    Over \(\ZZ_\ell\), the bound attributed by Serre to the method of Oesterl\'e and Robba is sharper. It is stated following the corollary to Theorem~8 in~\cite[\S3.2, p.~146]{serre}:
    Let \(F\in\ZZ_\ell[T_1,\ldots,T_N]\) be a nonzero polynomial of total degree \(D\), and define
    \begin{equation*}
        Y:=\{x\in\ZZ_\ell^N\mid F(x)=0\}\,.
    \end{equation*}
    For every integer \(m\ge 0\), let
    \begin{equation*}
        \r_m\colon\ZZ_\ell^N\to(\ZZ/\ell^m\ZZ)^N
    \end{equation*}
    be the coordinate-wise reduction map and define \(Y_m:=\r_m(Y)\).
    Then
    \[
    \#Y_m\le D\,\ell^{m(N-1)}
    \qquad\text{for every integer }m\ge 0.
    \]
\end{remark}

\begin{lemma}\label{lem:finitedyn}
    Let \(A,\pi,K,k,q,p\) be the same as in Lemma~\ref{lem:hypredcount}, and let \(\G\colon A^n\to A^n\) be of the form
    \[\G(x)=Jx+\pi B(x), \quad J\in \GL_n(A), \quad B\colon A^n\to A^n,\]
    where \(B\) preserves congruence modulo \(\pi^m\) for all \(m\ge 1\), that is,
    \[x\equiv y\pmod{\pi^m}\implies B(x)\equiv B(y) \pmod{\pi^m}.\]
    Let \(h\ge 1\) be the order of the reduction of \(J\) in \(\GL_n(k)\). Then for all \(m\ge 1\), \(\G\) induces a map \(\G_m\) on \((A/\pi^mA)^n\) satisfying
    \begin{equation}\label{perm}
        \G_m^{hp^{m-1}}=\id
    \end{equation}
    with the same \(h\) for every \(m\).
    That is, \(\G_m\) is a permutation on \((A/\pi^mA)^n\) with order \mbox{dividing \(hp^{m-1}\)}.
\end{lemma}
\begin{proof}
    \begin{samepage}
    To show that \(\G\) induces a map \(\G_m\) on \((A/\pi^mA)^n\), let \(x-y\in \pi^mA^n\).
    Then 
    \begin{align*}
        \G(x)-\G(y)&=J(x-y)+\pi(B(x)-B(y))\\
        &\in\pi^mA^n+\pi^{m+1}A^n\subseteq\pi^mA^n\,.
    \end{align*}
    Therefore \(\G_m\) is well-defined.
    \end{samepage}
    
    To prove equation~\eqref{perm}, we first show by induction on \(j\ge 1\) that if \(x-y\in \pi^mA^n\), then
    \begin{equation}\label{linear}
        \G^j(x)-\G^j(y)\equiv J^j(x-y)\pmod{\pi^{m+1}}.
    \end{equation}
    That is, \(\G\) behaves like the linear map \(J\) on classes of \(A^n/\pi^mA^n\) when reduced modulo \(\pi^{m+1}\).
    For the base case \(j=1\), since \(x-y\in \pi^m A^n\), we have
    \[\G(x)-\G(y)=J(x-y)+\pi (B(x)-B(y))\equiv J(x-y)\pmod{\pi^{m+1}},\]
    as \(B\) preserves congruence levels and \(B(x)-B(y)\in \pi^m A^n\).
    Now assume equation~\eqref{linear} holds for \(j-1\),
    which implies that
    \[\G^{j-1}(x)-\G^{j-1}(y)=J^{j-1}(x-y)+\pi^{m+1}e\]
    for some \(e\in A^n\).
    Then we compute
    \begin{equation*}
    \begin{aligned}
        \G^j(x)-\G^j(y)
        &=J(\G^{j-1}(x)-\G^{j-1}(y))\\
        &\quad+\pi(B(\G^{j-1}(x))-B(\G^{j-1}(y)))\\
        &=J^j(x-y)+\pi^{m+1}Je+\pi^{m+1}e'\,.
    \end{aligned}
    \end{equation*}
    Here \(e'\in A^n\) exists by the induction hypothesis and preservation of congruence levels by \(B\). Both error terms lie in \(\pi^{m+1}A^n\), proving~\eqref{linear}.

    Since \(J\in\GL_n(A)\), its reduction modulo \(\pi\) gives
    \(\overline{J}\in\GL_n(k)\), which has order \(h\) by definition.
    We will prove by induction on \(m\) that
    \[
        \G^{hp^{m-1}}(x)\equiv x\pmod{\pi^m}
        \qquad\text{for every }x\in A^n.
    \]
    When \(m=1\), we have \(\G(x)\equiv Jx\pmod{\pi}\), so
    \[
        \G^h(x)\equiv J^h x\equiv x\pmod{\pi},
    \]
    since \(\overline{J}^h=\id\).

    Now assume that the claim holds for some \(m\ge 1\), and let
    \(N=hp^{m-1}\).
    Fix \(x\in A^n\). By assumption we have
    \[
        \G^N(x)-x\in\pi^m A^n.
    \]
    For each integer \(j\ge 0\), we find
    \begin{align*}
        \G^{(j+1)N}(x)-\G^{jN}(x)
        &=\G^{jN}(\G^N(x))-\G^{jN}(x)\\
        &\equiv J^{jN}(\G^N(x)-x)\pmod{\pi^{m+1}}\\
        &\equiv \G^N(x)-x\pmod{\pi^{m+1}}.
    \end{align*}
    The first congruence is an equality when \(j=0\) and follows from
    equation~\eqref{linear} when \mbox{\(j\ge 1\)}.
    The second congruence follows because \(\G^N(x)-x\in\pi^m A^n\) and \(h\) divides \(jN\), so \mbox{\(J^{jN}\equiv\id\pmod{\pi}\)}.

    \begin{samepage}
    Summing these congruences for \(j=0,\ldots,p-1\), we find
    \begin{align*}
        \G^{pN}(x)-x
        &=\sum_{j=0}^{p-1}
        \bigl(\G^{(j+1)N}(x)-\G^{jN}(x)\bigr)\\*
        &\equiv p(\G^N(x)-x)\pmod{\pi^{m+1}}\,.
    \end{align*}
    \end{samepage}
    \begin{samepage}
    Since \(\characteristic(k)=p\), we have \(p\in\pi A\), and therefore
    \(p(\G^N(x)-x)\in\pi^{m+1}A^n\).
    Thus
    \[
        \G^{hp^m}(x)=\G^{pN}(x)\equiv x\pmod{\pi^{m+1}},
    \]
    finishing the induction.
    \end{samepage}

    Thus~\eqref{perm} holds for every \(m\ge 1\), so \(\G_m\) is a permutation
    of \((A/\pi^m A)^n\) whose order \mbox{divides \(hp^{m-1}\)}.
\end{proof}

\section{Proof of the main theorem}
\begin{mainrestatement}[restated]
   Let \(f\colon \PP^n\to \PP^n\) be a dominant endomorphism over a number field \(K\). Then there is a finite extension \(L/K\) such that for every proper closed subset \(V\subsetneq \PP^n\), there is a point \(P\in \PP^n(L)\) such that
    \begin{equation*}
        \grano_{f,L}(P)\cap V(L)=\varnothing\,.
    \end{equation*}
\end{mainrestatement}
\begin{proof}
    The assertion is immediate if \(n=0\), so assume that \(n\ge 1\). We first assume that \(\deg(f)\ge 2\).
    By Lemma~\ref{lem:algperpt}, there exists a repelling periodic point \(\a\in\PP^n(\Kbar)\) such that \(f^r(\a)=\a\) for some \(r\ge 1\). We replace \(f\) by \(f^r\) using Lemma~\ref{lem:iterinv} so that \(\a\) becomes a repelling fixed point of \(f\).
    Now take a finite extension \(L/K\) such that \(\a\in\PP^n(L)\) and \(i\in L\), where \(i\) is a square root of \(-1\). Choose coordinates over \(L\) so that \(\a=[1:0:\dots:0]\). In the affine patch \(X_0\ne 0\), let \(J=df_\a\). Since \(\a\) is repelling, \(\det(J)\ne 0\).

    Let \(V\subsetneq\PP^n\) be a proper closed subset, and take a hypersurface \(H\subsetneq\PP^n\) defined over \(L\) of degree \(d\ge 1\) such that \(V\subseteq H\).
    To find a grand orbit avoiding \(V\), it now suffices to find a grand orbit avoiding \(H\).

    Write \(f=[F_0:\dots:F_n]\), where \(F_i\) are homogeneous polynomials over \(L\) of degree \(\deg(f)\) with no common projective zeros over \(\overline L\).
    Choose a prime number \(p\) such that \(p\equiv3\pmod4\), \(p>d\), and at every prime of \(L\) above \(p\), the map \(f\) in these coordinates has good reduction and \(J\) is integral with unit determinant. These conditions exclude only finitely many primes from the infinitely many primes congruent to \(3\) modulo \(4\).
    Take \(\gp\mid p\) in \(\Ocal_L\), let \(A:=\Ocal_{L,\gp}\), and choose a uniformizer \(\pi\in A\).
    We use the valuation on \(L\) normalized so that \(\pi\) has valuation \(1\), and let \(k:=A/\pi A\).
    Since \(i\in L\) and \(p\equiv3\pmod4\), the polynomial \(T^2+1\) has a root in \(k\) but no root in \(\FF_p\), so \(\#k=q\ge p^2\).

    For any point \(x\in\PP^n(L)\), we may identify it with a point of \(\PP^n(A)\) by normalizing the coordinates so that their minimal valuation is \(0\).
    Good reduction makes the coordinate-wise reduction maps commute with \(f\), giving the following diagram for every \(m\ge 1\):
    \begin{equation}
    \begin{tikzcd}
    \PP^n(L) \arrow[r, "\cong"] \arrow[d, "f"] & \PP^n(A) \arrow[r, "\r_m"] \arrow[d, "f"] & \PP^n(A/\pi^mA) \arrow[d, "f_m"] \\
    \PP^n(L) \arrow[r, "\cong"]                & \PP^n(A) \arrow[r, "\r_m"]                & \PP^n(A/\pi^mA)
    \end{tikzcd}
    \end{equation}
    The map \(\r_m\) is surjective: a unit coordinate can be normalized to \(1\), and the remaining coordinates can be lifted to \(A\subseteq L\).
    Let \(H_m:=\r_m(H(L))\). Since any two points in the same grand \(f\)-orbit in \(\PP^n(L)\) reduce to points in the same grand \(f_m\)-orbit, it suffices to find a point in \(\PP^n(A/\pi^mA)\) whose grand \(f_m\)-orbit avoids \(H_m\).
    Each point of \(H(L)\) has a unit coordinate in its normalized coordinates. In each of the \(n+1\) corresponding affine charts, the equation of \(H\) is a nonzero polynomial of total degree at most \(d\). Applying Lemma~\ref{lem:hypredcount} in each chart gives
    \begin{equation}\label{projcount}
        \#H_m\le(n+1)(dq^{n-1})^m\,.
    \end{equation}

    We next construct a ball on whose reductions \(f\) acts as a permutation. Take the image \(a:=\r_1(\a)\) of the fixed point and consider the \(\pi\)-adic ball
    \begin{equation*}
        \Ucal:=\r_1^{-1}(a)\subseteq\PP^n(A)\,.
    \end{equation*}
    Note that \(\Ucal\cong A^n\) via the coordinate-wise identification
    \begin{align*}
        \iota\colon\Ucal&\to A^n\\
        [1:\pi x_1:\dots:\pi x_n]&\mapsto(x_1,\dots,x_n)\,.
    \end{align*}
    Since \(f\) fixes \(\a\), good reduction implies that \(f\) preserves \(\Ucal\). Conjugating by \(\iota\) therefore gives a well-defined map
    \begin{equation*}
        \G:=\iota\circ f\circ\iota^{-1}\colon A^n\to A^n\,.
    \end{equation*}
    We claim that \(\G\) is of the form
    \begin{equation*}
        \G(x)=Jx+\pi B(x),\qquad J\in\GL_n(A)\,,
    \end{equation*}
    where \(B\colon A^n\to A^n\) preserves congruence modulo \(\pi^m\) for all \(m\ge 1\).

    To see the claim, scale the homogeneous polynomials \(F_i\) so that their coefficients lie in \(A\) and their reductions have no common projective zeros. Since \(f(\a)=\a\), we have \(F_j(1,0)=0\) for \(1\le j\le n\) and \(F_0(1,0)\in A^\times\). Indeed, if \(F_0(1,0)\) were not a unit, all the reduced \(F_i\) would vanish at \(a\), contrary to good reduction. Scaling all \(F_i\) by this unit's inverse, we may assume that \(F_0(1,0)=1\).
    
    The affine coordinate functions of \(f\) are \(F_j(1,x)/F_0(1,x)\). Since \(F_0(1,0)=1\) and \mbox{\(F_j(1,0)=0\)}, the quotient rule shows that the linear terms of \((F_1(1,x),\dots,F_n(1,x))\) are \(Jx\). Substituting \(\pi x\) for \(x\) gives a factor of \(\pi\) in every nonconstant term and a factor of \(\pi^2\) in every term of degree at least two. Thus there are a polynomial \(C\in A[x_1,\dots,x_n]\) and a vector of polynomials \(D\in A[x_1,\dots,x_n]^n\) such that
    \begin{align*}
        F_0(1,\pi x)&=1+\pi C(x)\,,\\
        (F_1(1,\pi x),\dots,F_n(1,\pi x))&=\pi Jx+\pi^2D(x)\,.
    \end{align*}
    Since \(\iota\) divides the affine coordinates by \(\pi\), we obtain
    \begin{equation}\label{gammaform}
        \G(x)=\frac{Jx+\pi D(x)}{1+\pi C(x)}
        =Jx+\pi\frac{D(x)-C(x)Jx}{1+\pi C(x)}\,,
    \end{equation}
    where the last fraction defines \(B(x)\). The components of its numerator and the denominator are polynomials over \(A\), so their reductions modulo \(\pi^m\) depend only on the reduction of \(x\), as addition and multiplication commute with reduction. The denominator belongs to \(1+\pi A\) and is therefore a unit, whose inverse modulo \(\pi^m\) is also determined by its reduction. Thus \(B\) preserves congruence modulo every \(\pi^m\), proving the claim.

    By Lemma~\ref{lem:finitedyn}, there is an integer \(h\ge 1\), independent of \(m\), such that \(\G_m\) has order dividing \(hp^{m-1}\) for every \(m\ge 1\).
    Let \(\Ucal_m:=\r_m(\Ucal)\). For \(m\ge 2\), the coordinates above identify \(\Ucal_m\) with \((A/\pi^{m-1}A)^n\): reducing \(\pi x\) modulo \(\pi^m\) is equivalent to reducing \(x\) modulo \(\pi^{m-1}\). Thus
    \begin{equation*}
        \#\Ucal_m=q^{n(m-1)}\,.
    \end{equation*}
    The restriction \(f_m|_{\Ucal_m}\) corresponds to \(\G_{m-1}\) under this identification, so its order divides~\(hp^{m-2}\). We will use the weaker bound \(hp^{m-1}\).

    We must also account for points of \(H\) whose orbits enter the ball. We claim that any point that enters \(\Ucal\) under iteration of \(f\) does so within \(M-1\) steps, where \(M:=\#\PP^n(k)\).
    Let \(x\in\PP^n(L)\) and suppose that \(f^\ell(x)\in\Ucal\) for some \(\ell\ge 0\). By the commutative diagram,
    \begin{equation*}
        f_1^\ell(\r_1(x))=\r_1(f^\ell(x))=a\,.
    \end{equation*}
    Since \(a\) is fixed, the points preceding its first occurrence in this finite orbit are distinct. Otherwise, a repetition would trap the orbit in a cycle not containing \(a\). There are only \(M\) points in \(\PP^n(k)\), so \(f_1^{\ell'}(\r_1(x))=a\) for some \(0\le\ell'<M\). Thus \(f^{\ell'}(x)\in\Ucal\), proving the claim. The same argument after reduction shows that any point whose \(f_m\)-orbit enters \(\Ucal_m\) does so within \(M-1\) steps, since \(\Ucal_m\) consists exactly of the points reducing to \(a\) modulo \(\pi\).

    Write \(\grano_{f_m}(H_m)\) for the union of the grand \(f_m\)-orbits of points in \(H_m\). If \(x\in H_m\) and \(y\in\Ucal_m\) satisfy \(f_m^\ell(x)=f_m^j(y)\) for some \(\ell,j\ge 0\), then the orbit of \(x\) enters \(\Ucal_m\). Its first entrance occurs at some \(0\le\ell'<M\) with \(\ell'\le\ell\). Since \(f_m|_{\Ucal_m}\) is a permutation, its cycles are disjoint. The equality of iterates therefore places \(f_m^{\ell'}(x)\) and \(y\) on the same cycle in \(\Ucal_m\).
    Conversely, every point on a cycle in \(\Ucal_m\) meeting \(f_m^{\ell'}(H_m)\) belongs to \(\grano_{f_m}(H_m)\). Thus \(\grano_{f_m}(H_m)\cap\Ucal_m\) is precisely the union of the cycles in \(\Ucal_m\) meeting
    \begin{equation*}
        \bigcup_{j=0}^{M-1}\bigl(f_m^j(H_m)\cap\Ucal_m\bigr)\,.
    \end{equation*}

    It remains to show that these cycles do not exhaust \(\Ucal_m\). Since \(q\ge p^2>pd\), choose an integer \(m\ge 2\) such that
    \begin{equation*}
        m>\frac{\log\bigl(hM(n+1)q^n/p\bigr)}{\log\bigl(q/(pd)\bigr)}\,.
    \end{equation*}
    Then
    \begin{equation*}
        hp^{m-1}M(n+1)d^m q^n<q^m\,.
    \end{equation*}
    Each cycle in \(\Ucal_m\) has at most \(hp^{m-1}\) points. Since \(\#f_m^j(H_m)\le\#H_m\), equation~\eqref{projcount} gives
    \begin{align*}
        \#\bigl(\grano_{f_m}(H_m)\cap\Ucal_m\bigr)
        &\le hp^{m-1}\sum_{j=0}^{M-1}\#\bigl(f_m^j(H_m)\cap\Ucal_m\bigr)\\
        &\le hp^{m-1}M(n+1)(dq^{n-1})^m\\
        &<q^{n(m-1)}=\#\Ucal_m\,.
    \end{align*}
    Hence there exists a point \(x\in\Ucal_m\) whose grand \(f_m\)-orbit avoids \(H_m\). Lift its coordinates to a point \(P\in\Ucal\subseteq\PP^n(L)\). If \(\grano_{f,L}(P)\) met \(H(L)\), reducing the corresponding equality of iterates would contradict the choice of \(x\). Thus \(\grano_{f,L}(P)\cap V(L)=\varnothing\), finishing the proof for the \(\deg(f)\ge 2\) case.

    It remains to consider the case when \(\deg(f)=1\). Since \(f\) is dominant, there exists a matrix \(T\in\GL_{n+1}(K)\) such that \(f([x])=[Tx]\).
    Take \(L=K(i)\), where \(i^2=-1\), and let \(V\subsetneq\PP^n\) be a proper closed subset. As before, take a hypersurface \(H\) defined over \(L\) of degree \(d\ge 1\) containing \(V\). Choose a prime \(p\equiv3\pmod4\) with \(p>d\) such that \(T\) and \(T^{-1}\) are integral at every prime of \(L\) above \(p\). Take \(\gp\mid p\), let \(A:=\Ocal_{L,\gp}\) with uniformizer \(\pi\), and let \(k:=A/\pi A\) with \(\#k=q\ge p^2\).
    Since \(T\) and \(T^{-1}\) have entries in \(A\), they induce mutually inverse maps on \(\PP^n(A/\pi^mA)\) for every \(m\ge 1\).

    Let \(h\) be the order of the reduction of \(T\) in \(\GL_{n+1}(k)\). Thus there is a matrix \(R\in\Mat_{n+1}(A)\) such that \(T^h=\id+\pi R\).
    For every integer \(m\ge 1\), the binomial theorem gives
    \begin{equation}\label{pm}
        (T^h)^{p^{m-1}}=(\id+\pi R)^{p^{m-1}}
        =\id+\sum_{j=1}^{p^{m-1}}\binom{p^{m-1}}{j}(\pi R)^j\,.
    \end{equation}
    Write \(\nu_\pi\) for the valuation normalized by \(\nu_\pi(\pi)=1\), and \(\nu_p\) for the usual \(p\)-adic valuation on rational numbers. Since \(p\in\pi A\), for \(1\le j\le p^{m-1}\) we have
    \begin{equation}\label{coeff}
        \nu_\pi\left(\binom{p^{m-1}}{j}\pi^j\right)
        \ge m-1-\nu_p(j)+j\ge m\,.
    \end{equation}
    The first inequality follows from
    \begin{equation*}
        \binom{p^{m-1}}{j}=\frac{p^{m-1}}{j}\binom{p^{m-1}-1}{j-1}\,,
    \end{equation*}
    because the remaining binomial coefficient is an integer and \(\nu_\pi\ge\nu_p\) on nonzero integers. The second follows from \(j-\nu_p(j)\ge 1\).
    Thus \(T^{hp^{m-1}}\equiv\id\pmod{\pi^m}\). Consequently, the induced map \(f_m\) is a permutation of \(\PP^n(A/\pi^mA)\) with order dividing \(hp^{m-1}\).

    Let \(H_m:=\r_m(H(L))\). The same affine-chart argument as in~\eqref{projcount} gives
    \begin{equation*}
        \#H_m\le(n+1)(dq^{n-1})^m\,.
    \end{equation*}
    Since \(q\ge p^2>pd\), choose an integer \(m\ge 1\) such that
    \begin{equation*}
        m>\frac{\log\bigl(h(n+1)/p\bigr)}{\log\bigl(q/(pd)\bigr)}\,.
    \end{equation*}
    Then
    \begin{equation*}
        hp^{m-1}(n+1)d^m<q^m\,,
    \end{equation*}
    so, since every grand \(f_m\)-orbit is a cycle of length at most \(hp^{m-1}\),
    \begin{align*}
        \#\grano_{f_m}(H_m)
        &\le hp^{m-1}\#H_m\\
        &\le hp^{m-1}(n+1)(dq^{n-1})^m\\
        &<q^{nm}\le\#\PP^n(A/\pi^mA)\,,
    \end{align*}
    where the last inequality follows by counting the points in a single affine chart. Choose a point outside \(\grano_{f_m}(H_m)\) and lift its normalized coordinates to a point \(P\in\PP^n(L)\). As in the previous case, reduction of any equality of iterates shows that \(\grano_{f,L}(P)\cap H(L)=\varnothing\), and hence \(\grano_{f,L}(P)\cap V(L)=\varnothing\). This proves the assertion with the fixed extension \(L=K(i)\).
\end{proof}

\printbibliography

\end{document}